\documentclass[12pt,a4paper]{article}
\usepackage[utf8]{inputenc}
\usepackage[T1]{fontenc}
\usepackage{amsmath,amssymb,amsthm}
\usepackage{geometry}
\usepackage[colorlinks=true,linkcolor=blue,urlcolor=blue]{hyperref}

\newtheorem{theorem}{Theorem}
\newtheorem{basictheorem}{Basic Theorem}
\newtheorem{lemma}{Lemma}
\theoremstyle{remark}
\newtheorem{remark}{Remark}

\newcommand{\N}{\mathbb{N}}
\newcommand{\R}{\mathbb{R}}
\newcommand{\Z}{\mathbb{Z}}

\title{On The Extensions of the Binomial Theorem\\ And Properties Of The Binomial Coefficient}
\author{Shazali Abdalla Fadul\\
\small Faculty of Mathematical Sciences \& Statistics, Al-Neelain University, Khartoum, Sudan\\
\small e-mail: \texttt{shazlyabdullah3@gmail.com}}
\date{}

\begin{document}
\maketitle

\begin{abstract}
In this study we used an algebraic method that uses elementary algebra and the binomial theorem. Using this method, several results have been proven, which are considered an extension of the binomial theorem polynomial theorem
$L_n(x,y,a)=V_n^{\,n}(x,y,a)+S_n(x,y)$. It has several properties, for example if $a=1$ then $V_n^{\,n}(x,y,a)=0$, therefore $L_n(x,y,a)=S_n(x,y)$. From the results we obtained, if $a\in\R$ and $n$ is a positive integer, then
\[
\frac{(a^n-1)^2+2a^n-2-(a+1)^n(a^n-1)}{a}
= -\sum_{j=1}^{n-1}(a+1)^j\sum_{k=1}^{j}\binom{n}{k}(a-1)^k .
\]
In addition, we studied the properties of the binomial coefficient, in particular we studied the relationship between the coefficient of the binomial theorem and the Euler function $\phi(n)$ where $n\in\N$, proving the existence of an analytical relationship possessed by the binomial coefficient. If it is combined with the Euler function, from the results we prove that if $\gcd(2,n)=1$ and $\phi(n)$ is the Euler function, then
\[
\sum_{\substack{j\ge 1\\ 1\le j\le \phi(n)-1\\ j\ \text{is odd}}}\binom{\phi(n)}{j}3^{j}\equiv 0 \pmod{n2^{\phi(n)-1}} .
\]
Moreover, we obtained several results on prime numbers of the form $ka^m+1$; the work of H.\ C.\ Pocklington on these numbers has been generalized, with finite series. Also, in particular, the work of Zhi-Wei Sun and E.\ Mortenson on the properties of binomial coefficients has been generalized, as that work has been extended.
\end{abstract}


\section{Introduction}

The binomial theorem is one of the basic theorems in mathematics, where it plays a key role in proving a number of theorems in mathematical analysis, computation, and number theory. It also plays an essential role in proving many theorems in mathematical analysis.

\begin{theorem}[Binomial theorem]
Let $x,y\in\R$ and $n$ a positive integer. Then
\[
(x+y)^n=\sum_{j=0}^{n}\binom{n}{j}x^{j}y^{n-j}.
\]
\end{theorem}

\begin{theorem}[Difference of two $n$th powers]
Let $x,y\in\R$ and $n$ a positive integer (see \cite{rose}, p.\ 96). Then
\[
x^n-y^n=(x-y)\sum_{j=1}^{n}x^{n-j}y^{j-1}.
\]
\end{theorem}

In this study an algebraic method was used that uses elementary algebra, in particular the basic operations of primary algebra: multiplication, division, subtraction, addition, substitution, and extraction of the common factor, in extending a similar theorem, where we reached a theorem consisting of three variables that possesses many properties and features. We also concluded from this theorem, which is an extension of the binomial theorem, many theorems; for example, we proved that if $a\in\R$ and $n$ is a positive integer, then
\[
(a^{n-1}-1)(a^n-1)=\sum_{j=1}^{n-1}(a+1)^j\left(a^n-\sum_{k=1}^{j}\binom{n}{k}(a-1)^k\right).
\]

We also studied the properties of the binomial coefficient. In particular, we studied the relationship between the binomial coefficient and the Euler function, and we reached several conclusions. The study of arithmetic properties of binomial coefficients has a long history. In 1819 Babbage proved the congruence
\[
\binom{2p-1}{p-1}\equiv 1 \pmod{p^2}.
\]
That congruence is essentially an extension of Wilson's theorem: if $p$ is a prime then $(p-1)!\equiv -1 \pmod p$. This result was generalized by Gauss \cite{rose} to the following:
\[
\prod_{\substack{k=1\\ \gcd(k,m)=1}}^{m}k\equiv
\begin{cases}
-1 \pmod m & \text{if } m=4,\ p^{a},\ 2p^{a},\\
\ \ 1 \pmod m & \text{otherwise.}
\end{cases}
\]

In 2016 Guo-Shuai Mao and Zhi-Wei Sun \cite{maosun} proved that if $p>3$ is a prime and $a$ is a positive integer, then if $p\equiv 1 \pmod 4$ or $a>1$,
\[
\sum_{k=0}^{\lfloor \frac34 p^{a}\rfloor}\frac{\binom{2k}{k}^2}{16^k}\equiv\left(\frac{-1}{p^{a}}\right)\pmod{p^3}.
\]
And in 2003 Rodriguez-Villegas posed conjectures on
\[
\sum_{k=0}^{p-1}\frac{\binom{2k}{k}\binom{3k}{k}}{27^k}\equiv\left(\frac{p}{3}\right)\pmod{p^2},\qquad
\sum_{k=0}^{p-1}\frac{\binom{6k}{3k}\binom{3k}{k}}{432^k}\equiv\left(\frac{-1}{p}\right)\pmod{p^2}.
\]
These results were proven by E.\ Mortenson \cite{mort1}. In this study, this work was generalized to the following: if $n$ is a positive integer and $\phi(n)$ is the Euler function, then
\[
\sum_{\substack{j\ge 1\\ 1\le j\le \phi(m)-1\\ j\ \text{is odd}}}\binom{\phi(m)}{j}3^{j}\equiv 0\pmod{m2^{\phi(m)-1}}.
\]

The purpose of this study is first to study extensions of the binomial theorem; in particular, we used an algebraic method based on elementary algebra. The basic results are those that can be described as extensions of the binomial theorem. The properties of the binomial coefficient were also studied; in particular, we studied the relationship between the binomial coefficient and the Euler function $\phi(n)$. We reached several results on this topic, results that highlight that the relationship is given by polynomials with finite limits linking the binomial coefficient with the Euler function. Moreover, we studied the properties of prime numbers of the form $ka^m+1$; the work of Pocklington \cite{pock} was generalized for these numbers.

\section{Proof of the basic theorem}

In this section we study extensions of the binomial theorem. We use techniques that may be unexpected: the binomial theorem itself, elementary algebra, and the difference of two $n$th powers (see \cite{rose}, p.\ 96). We reached several results that differ from the binomial theorem in the number of variables as well as in the properties. In particular, the results are finitely bounded sums in three or fewer variables.

\begin{basictheorem}
Let $a,x,y$ be real numbers and $m$ a constant. Then
\[
L_n(x,y,a)=V_n^{\,n}(x,y,a)+S_n(x,y),
\]
where
\[
L_n(x,y,a)=(a-x+y)^n-m(1+y)^n,
\]
and
\[
V_n^{\,n}(x,y,a)=\sum_{j=0}^{n}(-1)^j\binom{n}{j}\left(a^{n-j}-m\right)(x-y)^j,
\]
and
\[
S_n(x,y)=mx\sum_{j=0}^{n-1}(-1)^{j+1}\frac{(x-y)^j}{y^{j+1}}\left((1+y)^n-\sum_{k=0}^{j}\binom{n}{k}y^k\right).
\]
\end{basictheorem}

\begin{proof}
Let $x,y$ be real numbers. According to the difference of two $n$th powers theorem we have
\[
(x-y)^n-(-y)^n=x\sum_{j=1}^{n}(x-y)^{j-1}(-y)^{n-j}.
\]
Then
\[
-(-y)^n=-(x-y)^n+x\sum_{j=1}^{n}(x-y)^{j-1}(-y)^{n-j}.
\]
Let $q\in\R$, $n\in\N$, where $m$ is a constant; then by multiplying by $m$ and adding $a^{q}(x-y)^n$ to both sides,
\begin{align*}
a^q(x-y)^n-m(-y)^n
&=a^q(x-y)^n-m(x-y)^n+mx\sum_{j=1}^{n}(x-y)^{j-1}(-y)^{n-j}\\
&=(a^q-m)(x-y)^n+mx\sum_{j=1}^{n}(x-y)^{j-1}(-y)^{n-j}. \tag{1}
\end{align*}
According to the binomial theorem,
\begin{align*}
(a-x+y)^n-m(1+y)^n
&=\sum_{j=0}^{n}(-1)^j\binom{n}{j}\left(a^{n-j}(x-y)^j-m(-y)^j\right). \tag{2}
\end{align*}
Now by substituting equation (1) into equation (2) we get
\begin{align*}
(a-x+y)^n-m(1+y)^n
&=\sum_{j=0}^{n}(-1)^j\binom{n}{j}\left((a^{n-j}-m)(x-y)^j+mx\sum_{k=1}^{j}(x-y)^{k-1}(-y)^{j-k}\right)\\
&=\sum_{j=0}^{n}\binom{n}{j}(a^{n-j}-m)(-x+y)^j
+mx\sum_{j=1}^{n}\sum_{k=1}^{j}\binom{n}{j}(-1)^j(x-y)^{j}(-y)^{j-k}. \tag{3}
\end{align*}
Now note that the new equation consists of two terms; the first term consists of ordered terms that cannot be arranged in another way, but the second term consists of double sums, which means that it is possible to arrange the terms in another way. Notice
\begin{align*}
mx\sum_{j=1}^{n}\sum_{k=1}^{j}\binom{n}{j}(-1)^j(x-y)^j(-y)^{j-k}
&=\sum_{j=1}^{n}(-1)^j(x-y)^{j-1}\frac{1}{y^j}\left(\sum_{k=j}^{n}\binom{n}{k}y^k\right)\\
&=\sum_{j=0}^{n}(-1)^j\frac{(x-y)^j}{y^{j+1}}\left((1+y)^n-\sum_{k=0}^{j}\binom{n}{k}y^k\right). \tag{4}
\end{align*}
Then from (3) and (4) we have
\begin{align*}
(a-x+y)^n-m(1+y)^n
&=\sum_{j=0}^{n}\binom{n}{j}(a^{n-j}-m)(-x+y)^j\\
&\quad+mx\sum_{j=0}^{n}(-1)^{j+1}\frac{(x-y)^j}{y^{j+1}}\left((1+y)^n-\sum_{k=0}^{j}\binom{n}{k}y^k\right).
\end{align*}
\end{proof}

\begin{remark}
Notice the proof process is a combinatorial process. The binomial theorem was merged with the difference of two $n$th powers theorem into one theorem. In the proof, we note first that we got equation (1) using the difference of two $n$th powers theorem, after that we reached equation (2) from the binomial theorem, then equation (1) was merged with (2), then the terms were arranged until we got the final equation. This differs from the proofs of the two theorems, which are proved using mathematical induction, see \cite{rose} and \cite{tatt}. If we make comparisons between the binomial theorem and Theorem~1, in fact we find that Theorem~1 is more general than the binomial theorem, in particular in terms of the number of variables and the properties that characterize it. For example, if we substitute $m=0$ into Theorem~1 with ease, we get the binomial theorem; this means the binomial theorem has become a special case of Theorem~1 when substituting certain values for the variable $m$. Moreover, note the difference of two $n$th powers theorem $x^n-y^n=(x-y)(x^{n-1}+x^{n-2}y+x^{n-3}y^2+\cdots+y^{n-1})$. Notice in Theorem~1 that $L_n(x,y,a)=V_n^{\,n}(x,y,a)+S_n(x,y)$ where $L_n(x,y,a)=(a-x+y)^n-m(y+1)^n$. Now, if we put these values $m=1$ and $a=1$ and $y=a-1$, we find that $L_n(x,a-1,a)=(x+a)^n-a^n$ and $V_n^{\,n}(x,a-1,0)=0$; then we have $(x+a)^n-a^n=S_n(x,a-1)$, and this means that this is a new analysis of the limits of $x^n-y^n$ other than the usual ones in the difference of two $n$th powers. Note $S_n(x,a-1)=x\left((x+a)^{n-1}+(x+a)^{n-2}a+\cdots+a^{n-1}\right)$, and this is one of the properties of Theorem~1. We notice that, because of the multiple variables in Theorem~1, it became rich in properties, because it is possible to deduce several formulas according to the values of those variables. We will explain these properties and advantages in proving the following theorems.
\end{remark}

\begin{theorem}
If $a,d>0$ are real numbers and $n\in\N$, then
\[
a^{n-1}(d^n-1)=(d-1)\sum_{j=0}^{n-1}\frac{(ad-1)^j}{(a-1)^{j+1}}\left(a^n-\sum_{k=0}^{j}\binom{n}{k}(a-1)^k\right).
\]
\end{theorem}

\begin{proof}
Let $a=1$ and $m=1$ in $V_n^{\,n}(x,y,a)$. Then we have
\[
V_n^{\,n}(x,y,1)=\sum_{j=1}^{n}(-1)^j\binom{n}{j}\left((1^{n-j}-1)(x-y)^j\right)=0.
\]
Then, according to the basic series,
\[
L_n(x,y,1)=V_n^{\,n}(x,y,1)+S_n(x,y),\qquad\text{so}\qquad L_n(x,y,1)=S_n(x,y).
\]
Then according to Basic Theorem~1 and from equation (4) we find that
\[
(1-x+y)^n-(1+y)^n=x\sum_{j=0}^{n-1}(-1)^{j+1}\frac{(x-y)^j}{y^{j+1}}\left((1+y)^n-\sum_{k=0}^{j}\binom{n}{k}y^k\right).
\]
Let $a,d\in\R$ where $y=a-1$ and $x=-ad+a$. Now some substitution of these values in the equation will cause fundamental changes in the values $y$ and $x$, which are algebraically related in the equation. Note that $(1-x+y)^n=(1+ad-a+a-1)^n=(ad)^n$ and $(x-y)=(-ad+a-a+1)=(-ad+1)$; then we have
\[
(ad)^n-a^n=(-ad+a)\sum_{j=0}^{n-1}(-1)^{j+1+j}\frac{(ad-1)^j}{(a-1)^{j+1}}\left(a^n-\sum_{k=0}^{j}\binom{n}{k}(a-1)^k\right).
\]
Now dividing by $a$ we have
\[
a^{n-1}(d^n-1)=(d-1)\sum_{j=0}^{n-1}\frac{(ad-1)^j}{(a-1)^{j+1}}\left(a^n-\sum_{k=0}^{j}\binom{n}{k}(a-1)^k\right).
\]
\end{proof}

\begin{remark}
Note that Theorem~2 was proved by substituting the values $a=1$, $m=1$ into $V_n^{\,n}(x,y,a)$ by the direct substitution method without performing any other operations. Remember this is one of the properties that characterizes Theorem~1, as we mentioned in the previous note, and this proves that, because of the number of variables the theorem has, it has become rich in properties. We note that analysis formulas were presented for the limits $a^{n-1}(d^n-1)$, which were deduced from Theorem~1 in an algebraic way. We note that the term $a^{n-1}$ is multiplied by the term $(d^n-1)$. If we transform the conclusion of this analysis from the binomial theorem $(a+d)^n$, we note that this is not possible as a comprehensive algebraic term $a+d$, and it is not possible to reduce the sequence of that term like the limit $a^{n-1}(d^n-1)$, which consists of $a^{n-1}$ and $(d^n-1)$. In addition, we note the analytical formulas for $a^{n-1}(d^n-1)$ differ from the analytical formulas for the binomial theorem in the formation of terms, but Theorem~2 participates in the field $\R$. For Theorem~2 the formulas are verified for each $\forall a,d\in\R$, and also the binomial theorem for all $a,d\in\R$.
\end{remark}

\begin{lemma}
Let $a\in\R$ and $n\in\N$. Then we have
\[
(a^{n-1}-1)(a^n-1)=\sum_{j=1}^{n-1}(a+1)^j\left(a^n-\sum_{k=0}^{j}\binom{n}{k}(a-1)^k\right).
\]
\end{lemma}

\begin{proof}
Let $a=d$ in Theorem~2, where $a>1$ is real. Then
\begin{align*}
a^{n-1}(a^n-1)
&=(a-1)\sum_{j=0}^{n-1}\frac{(a^2-1)^j}{(a-1)^{j+1}}\left(a^n-\sum_{k=0}^{j}\binom{n}{k}(a-1)^k\right)\\
&=(a-1)\sum_{j=0}^{n-1}\frac{\big((a-1)(a+1)\big)^j}{(a-1)^{j+1}}\left(a^n-\sum_{k=0}^{j}\binom{n}{k}(a-1)^k\right)\\
&=a^n-1-\sum_{j=1}^{n-1}(a+1)^j\left(a^n-\sum_{k=1}^{j}\binom{n}{k}(a-1)^k\right).
\end{align*}
Then
\[
(a^{n-1}-1)(a^n-1)=\sum_{j=1}^{n-1}(a+1)^j\left(a^n-\sum_{k=1}^{j}\binom{n}{k}(a-1)^k\right).
\]
\end{proof}

Lemma~1 is a special case of Theorem~2 when $a=d$. We notice some substitution of $a=d$ gave analytical formulas for the term $(a^{n-1}-1)(a^n-1)$. We notice that in Theorem~2 there is the term $ad-1$ but in Lemma~1 it contains the term $(a+1)$. Moreover, a boilerplate of the variable $d$ appears at the other end.

\begin{theorem}
If $a\in\R$ and $n\in\N$, then
\[
\frac{(a^n-1)^2+2a^n-2-(a+1)^n(a^n-1)}{a}=-\sum_{j=1}^{n-1}(a+1)^j\sum_{k=1}^{j}\binom{n}{k}(a-1)^k.
\]
\end{theorem}

\begin{proof}
According to Lemma~1 we have
\begin{align*}
(a^{n-1}-1)(a^n-1)&=\frac{(a^n-1)^2+a^n-1}{a}-(a^n-1)\\
&=\sum_{j=1}^{n-1}(a+1)^j(a^n-1)-\left(\sum_{j=1}^{n-1}(a+1)^j\sum_{k=1}^{j}\binom{n}{k}(a-1)^k\right)\\
&=\left(\frac{(a+1)^n-1}{a}-1\right)(a^n-1)-\left(\sum_{j=1}^{n-1}(a+1)^j\sum_{k=1}^{j}\binom{n}{k}(a-1)^k\right).
\end{align*}
Then we have
\begin{align*}
a^{2n-1}-a^{n-1}-a^n+1&=\frac{(a^n-1)^2+a^n-1}{a}-(a^n-1)\\
&=\left(\frac{(a+1)^n-1}{a}-1\right)(a^n-1)-\left(\sum_{j=1}^{n-1}(a+1)^j\sum_{k=1}^{j}\binom{n}{k}(a-1)^k\right).
\end{align*}
Then
\begin{align*}
&\frac{(a^n-1)^2+a^n-1}{a}-(a^n-1)-\left(\frac{(a+1)^n-1}{a}-1\right)(a^n-1)\\
&\qquad=\frac{(a^n-1)^2+2a^n-2-(a+1)^n(a^n-1)}{a}
=-\left(\sum_{j=1}^{n-1}(a+1)^j\sum_{k=1}^{j}\binom{n}{k}(a-1)^k\right).
\end{align*}
\end{proof}

\begin{theorem}
Let $a\in\R$ and $n\in\N$. Then
\begin{multline*}
\frac{(a^n-1)^2+2a^n-2-(a+1)^n(a^n-1)+\left((a+1)^n-1\right)\left(a^n-1-(a-1)^n\right)}{a}\\
=\sum_{j=1}^{n-1}\binom{n}{j}(a-1)^j\sum_{k=1}^{j}(a+1)^{k-1}.
\end{multline*}
\end{theorem}

\begin{proof}
From Theorem~3 we have
\begin{align*}
\frac{(a^n+1)^2+2a^n-2-(a+1)^n(a^n-1)}{a}
&=-\sum_{j=1}^{n-1}\binom{n}{j}(a-1)^j\left(\frac{(a+1)^n-1}{a}-\sum_{k=1}^{j}(a+1)^{k-1}\right)\\
&=-\sum_{j=1}^{n-1}\binom{n}{j}(a-1)^j\left(\frac{(a+1)^n-1}{a}\right)+\sum_{j=1}^{n-1}\binom{n}{j}(a-1)^j\sum_{k=1}^{j}(a+1)^{k-1}.
\end{align*}
Then we now reduce the term $\sum_{j=1}^{n-1}\binom{n}{j}(a-1)^j$ using the binomial theorem; note that $\sum_{j=1}^{n-1}\binom{n}{j}(a-1)^j=(a)^n-(a-1)^n-1$, so we get
\begin{multline*}
\frac{(a^n-1)^2-2a^n-2-(a+1)^n(a^n-1)}{a}+\frac{\left((a+1)^n-1\right)\left(a^n-1-(a-1)^n\right)}{a}\\
=\sum_{j=1}^{n-1}\binom{n}{j}(a-1)^j\sum_{k=1}^{j}(a+1)^{k-1}.
\end{multline*}
Then we have
\begin{multline*}
\frac{(a^n-1)^2+2a^n-2-(a+1)^n(a^n-1)+\left((a+1)^n-1\right)\left(a^n-1-(a-1)^n\right)}{a}\\
=\sum_{j=1}^{n-1}\binom{n}{j}(a-1)^j\sum_{k=1}^{j}(a+1)^{k-1}.
\end{multline*}
\end{proof}

\section{Congruence and prime numbers}

In this section we apply some of the theorems that we reached in the second section to the study of prime numbers, and they coincide in general, in particular in prime numbers. We study numbers of the form $ka^m+1$. We will generalize the work of Pocklington \cite{pock} about these numbers. In addition, we derive some theorems based on Euler's theorem.

\begin{theorem}
If $g\in\N$ where $g=hb^m+1$ and $\phi(g)$ is the Euler function and $1\le k\le\phi(g)$ where $k\mid\phi(g)$, then
\[
\left(hb^{m-k}\right)^{\frac{\phi(g)}{k}}\equiv 1\pmod g.
\]
\end{theorem}

\begin{proof}
From Theorem~2 we have that
\begin{align*}
a^{n-1}(d^n-1)-\frac{d-1}{a-1}(a^n-1)
&=(d-1)\sum_{j=1}^{n-1}\frac{(-1)^j(ad+1)^j}{(a-1)^{j+1}}\left(a^n-\sum_{k=0}^{j}\binom{n}{k}(a-1)^k\right).
\end{align*}
Note that
\[
(d-1)\sum_{j=1}^{n-1}\frac{(-1)^j(ad+1)^j}{(a-1)^{j+1}}\left(a^n-\sum_{k=0}^{j}\binom{n}{k}(a-1)^k\right)\equiv 0\pmod{ad+1}.
\]
Then
\[
a^{n-1}\left((-d)^n-1\right)-\frac{-d-1}{a-1}(a^n-1)\equiv 0\pmod{ad+1}.
\]
Therefore we have
\[
(-d)^n-1-(a^n-1)\equiv 0\pmod{ad+1}.
\]
Then
\[
(-d)^n\equiv 1\pmod{ad+1}\ \text{ if and only if }\ a^n\equiv 1\pmod{ad+1}.
\]
Now let $d=b^k$ and $a=hb^{m-k}$; then we have
\begin{equation}
(-b)^{kn}\equiv 1\pmod{ad+1}\ \text{ if and only if }\ \left(hb^{m-k}\right)^n\equiv 1\pmod{hb^m+1}. \tag{5}
\end{equation}
From Euler's theorem, see \cite{rose} and \cite{tatt}, we have: if $\gcd(n,a)=1$ and $\phi(n)$ is the Euler function then $a^{\phi(n)}\equiv 1\pmod n$. Then let, in equation (5), $g=hb^m+1$ and $\phi(g)=kn$ where $2\le k\le\phi(g)$; then we have in equation (5) that $(-b)^{\phi(g)}\equiv 1\pmod g$, but $\phi(g)$ is even, then $(-b)^{\phi(g)}\equiv 1\pmod g=b^{\phi(g)}\equiv 1\pmod g$. Then
\begin{equation}
b^{\phi(g)}\equiv 1\pmod g\ \text{ if and only if }\ \left(hb^{m-k}\right)^{\frac{\phi(g)}{k}}\equiv 1\pmod g. \tag{6}
\end{equation}
Then from Euler's theorem, see \cite{rose,tatt},
\[
\left(hb^{m-k}\right)^{\frac{\phi(g)}{k}}\equiv 1\pmod g.
\]
\end{proof}

We note that we relied on Euler's theorem to prove this theorem. Euler's theorem is one of the most important theorems in number theory, and it was basically a generalization of Fermat's little theorem. In 1610 Fermat wrote in a letter to Frenicle that whenever $p$ is prime, $p$ divides $a^{p-1}-1$ for all integers $a$ not divisible by $p$, a result now known as Fermat's little theorem. An equivalent formulation is the assertion that $p$ divides $a^p-a$ for all integers $a$, whenever $p$ is prime. The question naturally arose as to whether the primes are the only integers exceeding $1$ that satisfy this criterion. Euler went further than that. First, he proved the conjecture and also generalized it: if $n$ is a positive integer and $\gcd(a,n)=1$ where $\phi(n)$ is the Euler function, then $a^{\phi(n)}\equiv 1\pmod n$ (see \cite{tatt}); this is now known as Euler's theorem, where Fermat's little theorem became a special case of Euler's theorem when $n=p$ where $p$ is prime. For a related discussion see Kenneth H.\ Rosen \cite[p.\ 55]{rose}. Notice that in Theorem~5 we introduced new properties of Euler's theorem when $n=hb^m+1$.

\begin{theorem}
If $p=hb^m+1$ where $h<m$, then
\[
\left(hb^{m-h}\right)^{a^m}\equiv 1\pmod p\quad\text{if $p$ is a prime.}
\]
\end{theorem}

\begin{proof}
Let, in Theorem~5, $p=hb^m+1$ and $h<m$ and $h=k$, where $p$ is a prime. Then $\left(hb^{m-h}\right)^{a^m}\equiv 1\pmod g$.
\end{proof}

The prime numbers of the form $N=k2^n+1$ are called Proth numbers, and were studied by the French mathematician Fran\c{c}ois Proth. Proth proved: if $N=k2^n+1$ where $k$ is odd and $k<2^n$, if $a^{\frac{N-1}{2}}\equiv -1\pmod N$ for some $a\in\Z$ then $N$ is a prime; this is the property of the Proth numbers. Also the numbers of the form $N=Kp^n+1$ are called generalized Proth numbers, studied by Pocklington \cite{pock}; also see \cite{grau2} and \cite{grau1}. An important step was made in 1914 by Pocklington; his result is the first generalization of Proth's theorem suitable for numbers of the form $Kp^n+1$. Pocklington proved: if $N=Kp^n+1$ where $K$ is odd and $k<p^n$ where $p$ is a prime, if for some $a\in\Z$
\[
a^{N-1}\equiv -1\pmod N\quad\text{and}\quad \gcd\!\left(a^{\frac{N-1}{p}}-1,\,p\right)\ \text{[as in the original]},
\]
then $N$ is prime. We note that Pocklington studied the properties of prime numbers of the form $N=Kp^n+1$ where $k<p^n$ and $p$ is a prime has two conditions: the first $k<p^n$ and the second that $p$ is a prime. Notice that the prime numbers of the form $N=ka^m+1$ for every $k,a,m\in\N$, not including these conditions, have not been studied, or those results do not include them. In this study we studied the numbers of the form $N=ka^m+1$ for every $k,a,m\in\N$. Notice that we presented results about this topic in Theorem~6, which is more general than those previous studies.

\begin{theorem}
If $q^n=hb^m+1$ where $\phi(q^n)$ is the Euler function and $\phi(q^n)/k$ where $1\le k\le m$, then
\[
\left(hb^{m-k}\right)^{\frac{\phi(q^n)}{k}}\equiv 1\pmod{q^n}.
\]
\end{theorem}

\begin{proof}
Let, in Theorem~6, $g=q^n$.
\end{proof}

\begin{theorem}
If $m$ is a positive integer where $\phi(m)$ is the Euler function and $\gcd(a,m)=\gcd(a+1,m)=1$, then
\[
\sum_{j=1}^{\phi(m)-1}\binom{\phi(m)}{j}(a-1)^j\sum_{k=1}^{j}(a+1)^{k-1}\equiv 0\pmod m.
\]
\end{theorem}

\begin{proof}
Let, in Theorem~4, $n=\phi(m)$; then we have
\begin{multline*}
\frac{\left(a^{\phi(m)}-1\right)^2+2a^{\phi(m)}-2-(a+1)^{\phi(m)}\left(a^{\phi(m)}-1\right)+\left((a+1)^{\phi(m)}-1\right)\left(a^{\phi(m)}-1-(a-1)^{\phi(m)}\right)}{a}\\
=\sum_{j=1}^{\phi(m)-1}\binom{\phi(m)}{j}(a-1)^j\sum_{k=1}^{j}(a+1)^{k-1}.
\end{multline*}
From Euler's theorem, if $\gcd(a,m)=1$ then $a^{\phi(m)}\equiv 1\pmod m$; then if $\gcd(a,m)=\gcd(a+1,m)=1$ we have that
\[
\frac{\left(a^{\phi(m)}-1\right)^2+2a^{\phi(m)}-2-(a+1)^{\phi(m)}\left(a^{\phi(m)}-1\right)+\left((a+1)^{\phi(m)}-1\right)\left(a^{\phi(m)}-1-(a-1)^{\phi(m)}\right)}{a}\equiv 0\pmod m.
\]
Then
\[
\sum_{j=1}^{\phi(m)-1}\binom{\phi(m)}{j}(a-1)^j\sum_{k=1}^{j}(a+1)^{k-1}\equiv 0\pmod m.
\]
\end{proof}

\begin{theorem}
If $a\in\R$ where $n$ is a positive integer, then
\[
\sum_{j=1}^{n-1}(a+1)^j\sum_{k=1}^{j}\binom{n}{k}(a-1)^k\equiv 0\pmod{a^n-1}.
\]
\end{theorem}

\begin{proof}
From Theorem~3 we have that
\[
\frac{(a^n-1)^2+2a^n-2-(a+1)^n(a^n-1)}{a}=-\sum_{j=1}^{n-1}(a+1)^j\sum_{k=1}^{j}\binom{n}{k}(a-1)^k.
\]
Then
\[
\frac{(a^n-1)^2+2a^n-2-(a+1)^n(a^n-1)}{a}\equiv 0\pmod{a^n-1}.
\]
Therefore
\[
\sum_{j=1}^{n-1}(a+1)^j\sum_{k=1}^{j}\binom{n}{k}(a-1)^k\equiv 0\pmod{a^n-1}.
\]
\end{proof}

\section{Properties of the binomial coefficient}

In this section the properties of the binomial coefficient were studied; in particular we studied the relationship between the Euler function and the binomial coefficient. We found an analytical relationship between the binomial coefficient and the Euler function. We proved this relationship in the form of a polynomial that combines the Euler function and the binomial coefficient in one theorem. It is also possible to describe these results as a natural generalization of the results that were reached in Theorems~5 and~7 in particular. These results were generalized in a radical way by linking the binomial coefficient with the Euler function.

\begin{theorem}
If $n>1$, $n\in\N$, where $n$ is even, then
\[
2^{n-1}(2^n-1)=\sum_{\substack{j\ge 1\\ 1\le j\le n-1\\ j\ \text{is odd}}}\binom{n}{j}3^j.
\]
\end{theorem}

\begin{proof}
Note that from Theorem~3, if $a=-2$ then we have
\[
\frac{(2^n-1)^2+2(2^n-1)-(2^n-1)}{-2}=\sum_{j=1}^{n-1}\binom{n}{j}(-2-1)^j\sum_{k=1}^{j}(-2+1)^{k-1}.
\]
Then
\begin{equation}
\frac{(2^n-1)^2+(2^n-1)}{-2}=-2^{n-1}(2^n-1)=\sum_{j=1}^{n-1}\binom{n}{j}(-3)^j\sum_{k=1}^{j}(-1)^{k-1}. \tag{7}
\end{equation}
Note that
\begin{equation}
\sum_{k=1}^{j}(-1)^{k-1}=
\begin{cases}
0 & \text{if $j$ is even},\\
1 & \text{if $j$ is odd}.
\end{cases} \tag{8}
\end{equation}
Then from (7) and (8) we have
\[
2^{n-1}(2^n-1)=\sum_{\substack{j\ge 1\\ 1\le j\le n-1\\ j\ \text{is odd}}}\binom{n}{j}3^j.
\]
\end{proof}

\begin{theorem}
If $m\in\N$ where $m$ is odd and $\phi(m)$ is the Euler function, then
\[
\sum_{\substack{j\ge 1\\ 1\le j\le \phi(m)-1\\ j\ \text{is odd}}}\binom{\phi(m)}{j}3^{j-1}\equiv 0\pmod{m2^{\phi(m)-1}}.
\]
\end{theorem}

\begin{proof}
Let $n=\phi(m)$ in Theorem~10, so we get
\[
2^{\phi(m)-1}\left(2^{\phi(m)}-1\right)=\sum_{\substack{j\ge 1\\ 1\le j\le \phi(m)-1\\ \gcd(j,2)=1}}\binom{\phi(m)}{j}3^j.
\]
From Euler's theorem, see \cite{tatt}, we have that $2^{\phi(m)}\equiv 1\pmod m$ if $\gcd(2,m)=1$, so we get that
\[
2^{\phi(m)-1}\left(2^{\phi(m)}-1\right)\equiv 0\pmod{m2^{\phi(m)-1}}.
\]
Then
\[
\sum_{\substack{j\ge 1\\ 1\le j\le \phi(m)-1\\ \gcd(j,2)=1}}\binom{\phi(m)}{j}3^j\equiv 0\pmod{m2^{\phi(m)-1}}.
\]
\end{proof}

The study of the binomial coefficient has a long history, and the reason for that is the properties and characteristics of several binomial coefficients rich in properties. For example, Jacobi proved that
\[
\binom{2(p-1)/2}{(p-1)/3}\equiv -r\pmod p
\]
where $p\equiv 1\pmod 6$, $4p=r^2+27s^2$ and $r\equiv 1\pmod 3$. In 1852 Wolstenholme proved his famous theorem which states that if $p\ge 5$ is a prime then $\sum_{k=1}^{p-1}\frac1k\equiv 0\pmod{p^2}$. Another form of Wolstenholme's theorem, which can be easily deduced from Wolstenholme's theorem, is that $\binom{2p-1}{p-1}\equiv 1\pmod{p^3}$. There are several studies on this subject, see for example \cite{guo,maosun,mort1}. But we notice that all these studies agree about the presence of the number in the form $p^k$ where $p$ is a prime and its divisibility, and $n$ a positive integer, $n\in\N$, has not been studied or those studies do not include it. We note in the theorem that we studied the divisibility of $n$, all $n\in\N$; note we proved that
\[
\sum_{\substack{j\ge 1\\ 1\le j\le \phi(m)-1\\ j\ \text{is odd}}}\binom{m}{j}3^{j-1}\equiv 0\pmod{m2^{\phi(m)-1}}
\]
for all $n\in\N$. This theorem is generalized to the work of Guo-Shuai Mao and Zhi-Wei Sun \cite{maosun}; note that in this theorem we found a relationship between the binomial coefficient and the Euler function.

\begin{lemma}
Let $p$ be a prime number. Then
\[
\sum_{\substack{j\ge 1\\ 1\le j\le p-2\\ j\ \text{is odd}}}\binom{p-1}{j}3^{j-1}\equiv 0\pmod{p2^{p-2}}.
\]
\end{lemma}

\begin{proof}
Let, in Theorem~11, $n=p$ where $p$ is a prime.
\end{proof}

\begin{theorem}
If $n$ is an odd positive integer and $\phi(n)$ is the Euler function, then
\[
\sum_{\substack{j\ge 2\\ 2\le j\le \phi(n)\\ j\ \text{is even}}}\binom{\phi(n)}{j}3^j\equiv 0\pmod n.
\]
\end{theorem}

\begin{proof}
Note we have from Euler's theorem $4^{\phi(n)}\equiv 1\pmod n$; then from Theorem~11
\[
\sum_{\substack{j\ge 1\\ 1\le j\le \phi(n)-1\\ j\ \text{is odd}}}\binom{\phi(n)}{j}3^j-\left(4^{\phi(n)}-1\right)\equiv 0\pmod n.
\]
Now note that
\[
\sum_{\substack{j\ge 1\\ 1\le j\le \phi(n)-1\\ \gcd(j,2)=1}}\binom{\phi(n)}{j}3^j-\sum_{j=1}^{\phi(n)}\binom{\phi(n)}{j}3^j\equiv 0\pmod n.
\]
Then
\[
\sum_{\substack{j\ge 2\\ 2\le j\le \phi(n)\\ j\ \text{is even}}}\binom{\phi(n)}{j}3^j\equiv 0\pmod n.
\]
\end{proof}

\section*{Conclusion}

We used a certain algebraic method in proving the basic Theorem~1. After that, we presented analysis formulas for the term $a^{n-1}(d^n-1)$ as well as analysis formulas for the term $(a^{n-1}-1)(a^n-1)$. The formulas are a finite set of terms. These terms cannot be deduced by analysis formulas using the binomial theorem. In addition, the results that we reached from the basic Theorem~4 were applied to prime numbers. Primality tests were presented for the prime numbers of the form $ka^m+1$, and this result is considered a generalization of the work of H.\ C.\ Pocklington. We also studied the properties of analytical congruence, which we note at the end of the second section. In addition to that, we studied the properties of the binomial coefficient, and the work of Mao, Sun and E.\ Mortenson was expanded, where we studied the relationship between the binomial coefficient and the Euler function, and the properties of a particular type of prime number.

\end{document}